\documentclass[preprint]{elsarticle}

\usepackage[utf8]{inputenc}  
\usepackage[T1]{fontenc}     
\usepackage{lmodern}         
\usepackage{microtype}       
\usepackage[scaled]{helvet} 
\usepackage[english]{babel}  
\usepackage{graphicx}        
\usepackage{float}
\usepackage{amsthm}
\usepackage[a4paper, marginparsep=3mm, marginparwidth=10mm]{geometry}
\usepackage[bookmarks=false]{hyperref}            

\usepackage{tikz}
\usepackage{tikz-3dplot}
\usepackage{makecell}

\newtheorem{definition}{Definition}
\newtheorem{example}{Example}
\usepackage[bookmarks=false]{hyperref}            
\usepackage{comment}         
\biboptions{sort&compress}   

\journal{TISSEC}

\usepackage{amsmath,amsfonts,amssymb}

\usepackage{subcaption}

\usepackage{xcolor}

\graphicspath{{fig/}}

\begin{document}

\begin{frontmatter}

\title{Persistent Homology: A Pedagogical Introduction with Biological Applications}

\author[label1,label2]{Aurelie Jodelle Kemme\corref{cor1}}
\author[label1]{Collins A. Agyingi\corref{cor2}}

\address[label1]{Department of Mathematical Sciences, University of South Africa, Pretoria, South Africa}
\address[label2]{African Institute for Mathematical Sciences, Research and Innovation Centre, Kigali, Rwanda}

\cortext[cor1]{Emails: akemme@aimsric.org/jkemme@quantumleapafrica.org}
\cortext[cor2]{Email: agyinca@unisa.ac.za}
\begin{abstract}

Persistent Homology (PH) is a fundamental tool in computational topology, designed to uncover the intrinsic geometric and topological features of data across multiple scales. Originating within the broader framework of Topological Data Analysis (TDA), PH has found diverse applications ranging from protein structure and knot analysis to financial domains such as Bitcoin behavior and stock market dynamics. Despite its growing relevance, there remains a lack of accessible resources that bridge the gap between theoretical foundations and practical implementation for beginners. This paper offers a clear and comprehensive introduction to persistent homology, guiding readers from core concepts to real-world application. Specifically, we illustrate the methodology through the analysis of a 3-1 supercoiled DNA structure. The paper is tailored for readers without prior exposure to algebraic topology, aiming to demystify persistent homology and foster its broader adoption in data analysis tasks.
\end{abstract}
\end{frontmatter}


\section{Introduction}

Persistent Homology (PH) is a powerful topological tool within the field of Topological Data Analysis (TDA), emerging from foundational works in applied algebraic topology and computational geometry~\cite{edelsbrunner2000topological, zomorodian2005computing}. Developed to address limitations in traditional machine learning (ML) methods\cite{pun2018persistenthomologybasedmachinelearningapplications,carlsson2012topological}, PH offers a means to capture the global geometric and topological structure of complex data by identifying and tracking features such as loops, cycles, and voids across multiple scales~\cite{carlsson2009topology,aggarwal2023tight,reani2022cycle}.
Conventional Machine Learning (ML) approaches—whether supervised, unsupervised, or reinforcement learning—primarily rely on statistical relationships, distances, or local density estimations within datasets. Supervised learning focuses on leveraging labelled data to predict outcomes \cite{muhammad2015supervised,nasteski2017overview,farrelly2023shape}, unsupervised learning seeks patterns and groupings without predefined labels \cite{hahne2008unsupervised,naeem2023unsupervised}, and reinforcement learning trains agents to interact with environments to maximize rewards~\cite{goodfellow2016deep,kaelbling1996reinforcement}. While effective in many scenarios, these paradigms often struggle to reveal the deeper, intrinsic structure of high-dimensional or non-linear data. Persistent Homology complements unsupervised learning by shifting the focus from local patterns to the global topological invariants of data—properties that remain unchanged under continuous transformations~\cite{ghrist2008barcodes}. Unlike standard clustering techniques, which typically depend on Euclidean distances or local density measures~\cite{ester1996density}, PH constructs a sequence of nested simplicial complexes known as a filtration. Through this process, PH systematically captures the birth and death of topological features across scales, providing insights into the underlying shape and connectivity of data that traditional methods overlook. For instance, in clustering tasks, conventional algorithms may effectively group points based on proximity but fail to detect essential structures such as holes or voids, which could signify meaningful gaps or cycles within the data. Persistent Homology, by contrast, identifies these features by analyzing how connected components, loops, and higher-dimensional voids evolve as the scale parameter changes~\cite{otter2017roadmap}. Visualization tools such as persistence barcodes and persistence diagrams further enhance interpretability, offering a clear representation of the lifespan of these topological features~\cite{edelsbrunner2010computational}. While PH begins with pairwise Euclidean distances to construct simplicial complexes (e.g., Vietoris–Rips complexes), its strength lies in transcending simple metric-based approaches, revealing multi-scale topological signatures that enrich data understanding far beyond what conventional ML techniques can provide. In this context, Persistent Homology stands out as an innovative and robust method, bridging gaps in classical data analysis by incorporating topological perspectives into the exploration of complex datasets.

\subsection{Overview of the paper}

This paper offers an accessible introduction to persistent homology, a central concept in topological data analysis, using intuitive examples built from point cloud data. Our goal is to guide the reader step by step from basic topological ideas to their powerful application in data analysis. We begin on the second section by laying the mathematical foundation, introducing key notions such as topological spaces, homeomorphisms, homotopy, and metric spaces that allow us to understand shapes beyond rigid geometry. This section also introduces the idea of a point cloud—a collection of data points in space—and the role of distance functions in measuring similarity between them. This sets the stage for constructing meaningful topological structures from raw data. In the Third section, we explore simplicial complexes and homology. Simplicial complexes form the building blocks used to approximate shapes in data. These combinatorial structures enable us to capture the geometry and connectivity of the data in a form that can be analyzed algebraically. Next, but still in the same section, we delve into chain complexes and homology, which allow us to detect topological features. This step transforms our geometric constructions into algebraic objects that can be computed. The fourth section introduces the core topic of the paper: persistent homology. Here, we show how topological features evolve as we vary the scale at which we view the data. Using techniques like filtration, barcodes, and persistence diagrams, we learn how to capture and visualize multi-scale structures in complex datasets. The fifth section presents a real-world implementation of persistent homology applied to the 3-1 supercoiled DNA structure. Finally, the last section presents the conclusion and some future research directions. This paper is intended for readers with little to no background in algebraic topology. Our focus is on clarity and intuition, supported by real-world examples to demonstrate how persistent homology can uncover hidden structure in data.


\section{Mathematical Preliminaries}
We familiarize the reader with mathematical concepts from topology and metric spaces, forming the backbone of persistent homology.
\subsection{Topology}

Topology is a branch of mathematics that focuses on the properties of shapes that stay the same when they are stretched, bent, or twisted—but not torn or glued. This means that two shapes are considered the same in topology if one can be smoothly transformed into the other without cutting or attaching new parts.

A classic example is the comparison between a coffee mug and a doughnut: although they look very different, both have exactly one hole. In the eyes of topology, this makes them equivalent because their essential structure—what mathematicians call the {\it topological genus}—is the same.

\begin{figure}[H]
    \centering
    \includegraphics[width=0.9\linewidth]{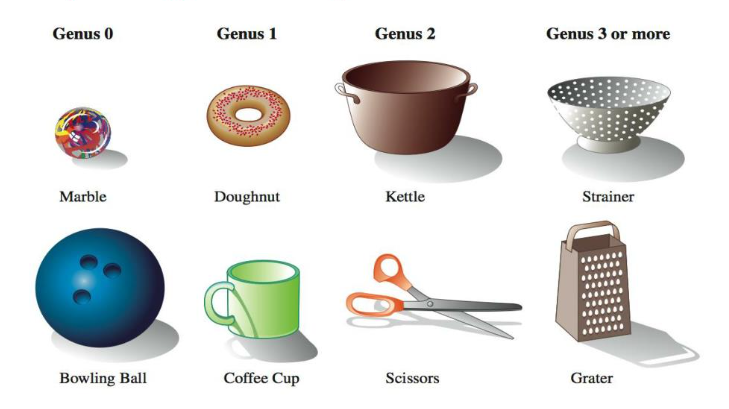}
    \caption{Some topologically equivalent objects. source: \href{https://www.slideserve.com/reese-cabrera/section-9-6-topology-powerpoint-ppt-presentation}{Reese Cabrera}}
    \label{coffemug}
\end{figure}

The image in Figure \ref{coffemug} offers a fun and intuitive way to understand the concept of {\bf genus} in topology. The {\bf genus} of a shape is a count of how many "holes" or "handles" it has—more formally, it's the number of independent loops you can draw on the object without being able to shrink them down to a point.

The categories in the image are:

\begin{enumerate}
    \item \textbf{Genus 0}: These objects have \textbf{no holes}. For example, a \textbf{marble} or a \textbf{bowling ball} is topologically equivalent to a sphere. You can’t pass anything through them—they are "solid."

    \item \textbf{Genus 1}: These shapes have \textbf{one hole}. The classic example is a \textbf{doughnut} (or torus), but interestingly, a \textbf{coffee cup} also fits here. Why? Because the handle of the cup forms a loop that is topologically the same as the hole in a doughnut.

    \item \textbf{Genus 2}: Objects like \textbf{kettles} and \textbf{scissors} have \textbf{two holes}—think of them as having two handles or loops.

    \item \textbf{Genus 3 or more}: Shapes such as a \textbf{strainer} or \textbf{grater} have \textbf{three or more holes}. These are more complex in their structure and can’t be continuously deformed into simpler ones without "breaking" them.
\end{enumerate}

This classification helps topologists understand which shapes are "the same" from a topological point of view, even if they look quite different geometrically. It's a key idea behind why a mug and a doughnut are topological twins!

Before we can understand how topology helps us study data using tools like persistent homology, we need to take a small step into the mathematical language of topology. Just as geometry is built on points, lines, and angles, topology is built on the idea of a \textbf{topological space}. A topological space provides a flexible framework to describe the notions of \textit{closeness} and \textit{connectedness} without relying on precise measurements like distance. This allows us to study the \textit{shape} or \textit{structure} of objects in a very general way, focusing on properties that remain unchanged when we bend or stretch them. Let us now introduce this fundamental concept more formally:

\begin{definition}[Topological Space]
A \textbf{topological space} $(X, \tau)$ is a set $X$ together with a collection of subsets $\tau$, called a \textit{topology}, satisfying the following properties:
\begin{enumerate}
    \item $\emptyset \in \tau$ and $X \in \tau$;
    \item $\tau$ is closed under finite intersections: if $A, B \in \tau$, then $A \cap B \in \tau$;
    \item $\tau$ is closed under arbitrary unions: if $\{A_i\}_{i \in I} \subseteq \tau$, then $\bigcup_{i \in I} A_i \in \tau$.
\end{enumerate}
The sets in $\tau$ are called \textit{open sets}.
\end{definition}

Once we have the notion of a \textit{topological space}, we can begin to compare different shapes by asking whether they are \textit{the same} from a topological point of view. But what does "the same" really mean in topology?

One of the strongest ways to consider two spaces as equivalent is through a concept called a \textbf{homeomorphism}, which is a continuous deformation between two spaces without cutting or gluing. 

\begin{definition}[Homeomorphism]
Let $X$ and $Y$ be two topological spaces. A function $f: X \rightarrow Y$ is called a \textbf{homeomorphism} if:
\begin{enumerate}
    \item $f$ is \textit{bijective} (one-to-one and onto),
    \item $f$ is \textit{continuous},
    \item The inverse function $f^{-1}: Y \rightarrow X$ is also \textit{continuous}.
\end{enumerate}
\end{definition}

If such a deformation exists, the two spaces are said to be \textit{homeomorphic}—for example, the famous case where a coffee mug is homeomorphic to a doughnut because both have one hole. However, in many situations, this notion can be too strict. This is where the idea of \textbf{homotopy} comes into play. 
In topology and algebra, \textbf{homotopy} provides a more flexible, or \textit{weaker}, notion of equivalence than homeomorphism. It allows us to classify shapes that may not be exactly deformable into each other but still share the same fundamental \textit{topological features}.

\begin{definition}[Homotopy]
Let $X$ and $Y$ be topological spaces, and let $f, g: X \rightarrow Y$ be two continuous functions. A \textbf{homotopy} between $f$ and $g$ is a continuous function:
\[
H: X \times [0,1] \rightarrow Y
\]
such that:
\[
H(x, 0) = f(x) \quad \text{and} \quad H(x, 1) = g(x) \quad \text{for all } x \in X.
\]
If such a function $H$ exists, we say that $f$ and $g$ are \textbf{homotopic}, denoted by $f \simeq g$.
\end{definition}

When two spaces can be continuously transformed into each other through a family of continuous maps, we say they are \textit{homotopically equivalent}. Homotopy focuses less on precise shape and more on preserving key features like connectedness and the number of holes. For example, a solid disk and a single point are homotopically equivalent because the disk can be continuously "shrunk" down to a point without tearing it.

\begin{definition}[Homotopy Equivalence]
Two topological spaces $X$ and $Y$ are said to be \textbf{homotopy equivalent} if there exist continuous functions:
\[
f: X \rightarrow Y \quad \text{and} \quad g: Y \rightarrow X
\]
such that:
\[
g \circ f \simeq \text{id}_X \quad \text{and} \quad f \circ g \simeq \text{id}_Y,
\]
where $\text{id}_X$ and $\text{id}_Y$ are the identity maps on $X$ and $Y$, respectively.

In this case, $X$ and $Y$ are considered to have the same "shape" from the perspective of homotopy theory, even if they are not homeomorphic.
\end{definition}

This concept of homotopy is essential in understanding how topology abstracts away geometric details to focus on the underlying structure, which will be crucial as we move toward tools like \textbf{homology} and \textbf{persistent homology} in data analysis.

\begin{figure}[H]
    \centering
    \includegraphics[width=0.5\linewidth]{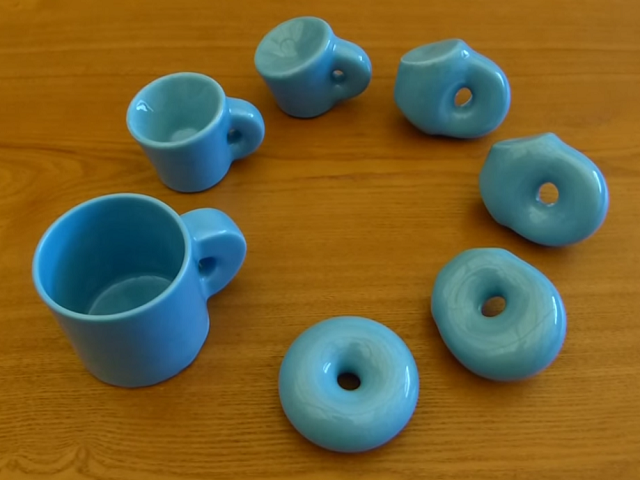}
     \caption{A coffee mug and a donut are homotopically equivalent. source: \href{https://www.gmanetwork.com/news/scitech/science/583886/when-is-a-coffee-mug-a-donut-topology-explains-it/story/\#goog_rewarded}{Marlowe Hood, GMA News}}

    \label{coffee mug}
\end{figure}

In real-world situations, data is rarely simple or neatly organized. Instead, it often forms complex patterns that are difficult to interpret using traditional methods. Understanding the overall structure—or \textit{shape}—of such data can quickly become challenging. This is exactly where the field of \textbf{Topological Data Analysis} (TDA) proves invaluable.

TDA, a growing area within data science and artificial intelligence, emerged from advances in applied algebraic topology and computational geometry. Its main purpose is to study the \textit{shape of data}, based on the idea that data isn't just a collection of points—it has structure, and this structure carries meaningful information. But how can we extract insights from the shape of data? To answer this, TDA offers several powerful tools, with two of the most well-known being the \textbf{Mapper algorithm} and \textbf{Persistent Homology (PH)}. The \textbf{Mapper algorithm} works by covering the dataset and connecting overlapping regions to form a simplified representation called a \textit{simplicial complex}. Think of it as creating a network or skeleton that captures the main features of the data's shape. However, using Mapper requires careful choices, like selecting the right function to project the data and deciding how to cover it, which can influence the results significantly. \textbf{Persistent Homology}, on the other hand, takes a different and often more systematic approach. Instead of relying on user-defined functions and covers, PH starts with data represented as a \textbf{point cloud}—a set of points in space. It calculates the distances between these points (usually using Euclidean distance) and builds a series of increasingly connected shapes called \textit{simplicial complexes}. This process, known as a \textbf{filtration}, allows us to observe how topological features like \textbf{connected components}, \textbf{loops}, and \textbf{voids} appear and disappear as we "zoom out" on the data.
The beauty of this method lies in its ability to detect which features are significant and which are just noise. Thanks to concepts from algebraic topology—particularly \textbf{homology}—we can track these features across different scales. The number of features in each dimension is counted using \textbf{Betti numbers}, giving us a clear summary of the data's shape. Finally, to make these insights more accessible, Persistent Homology provides visual tools such as \textbf{persistence barcodes} and \textbf{persistence diagrams}. These visuals show how long each topological feature "persists" throughout the filtration, helping us distinguish meaningful patterns from short-lived noise. It's important to note that the features we study belong to mathematical objects called \textbf{homology groups}, which, in this context, can be thought of as vector spaces where each vector represents a topological feature.

Before diving deeper into how Persistent Homology works, we must introduce two essential concepts that form the foundation of this method: the idea of a \textbf{point cloud} and the role of \textbf{metric spaces}. These will help us understand how to build a successful persistent homology pipeline.

\subsection{Metric Spaces}
Metric Spaces are a fundamental ideas in mathematics that formalizes what we intuitively think of as "distance". Before formally defining metric space and explaining its role in persistent homology computation, we first need to understand one fundamental concept: \textbf{point clouds}. Points cloud serves in general as the starting state for many techniques in Topological Data Analysis, including Persistent Homology.

\subsubsection{Point Clouds}

A \textbf{point cloud} is simply a collection of individual points scattered in space. You can think of it as a digital snapshot of an object or shape, where instead of smooth surfaces or lines, we only have dots representing key positions. Point clouds are often used to capture the structure of objects in three dimensions (3D), but they can exist in any number of dimensions depending on the data. Typically, each point in a point cloud is defined by coordinates—for example, $(X, Y, Z)$ in a 3D space. These points might come from various sources: they could be generated by 3D scanners, simulations, or even randomly, depending on whether the underlying distribution is known or unknown. While point clouds give us a way to visualize shapes or patterns, it's important to note that they don't naturally include connections or order between points. They are just scattered data, which means they lack an inherent structure. This is where mathematics helps us make sense of them. Point clouds can take many forms, depending on how the data is distributed. Some common examples include:

\begin{itemize}
    \item \textbf{Random Point Cloud}: Points are scattered without any clear pattern, often generated randomly.
    \item \textbf{Clusters}: Points group together in smaller regions, revealing natural groupings or categories within the data.
    \item \textbf{Spiral or Ellipsoid Shapes}: Points are arranged in a curved or circular pattern, suggesting more complex structures.
\end{itemize}

\begin{figure}[H]
    \centering
    \begin{subfigure}{0.32\textwidth}
        \centering
        \includegraphics[width=\linewidth]{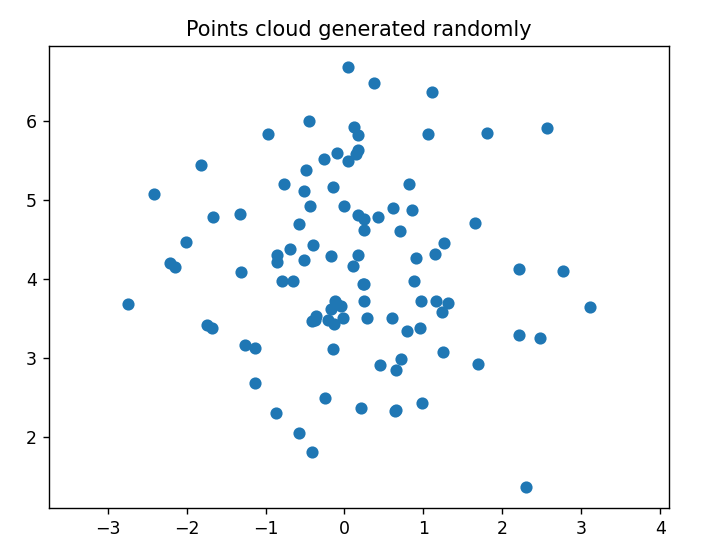} 
        \caption{Random points cloud}
        \label{fig:random}
    \end{subfigure}
    \hfill
    \begin{subfigure}{0.32\textwidth}
        \centering
        \includegraphics[width=\linewidth]{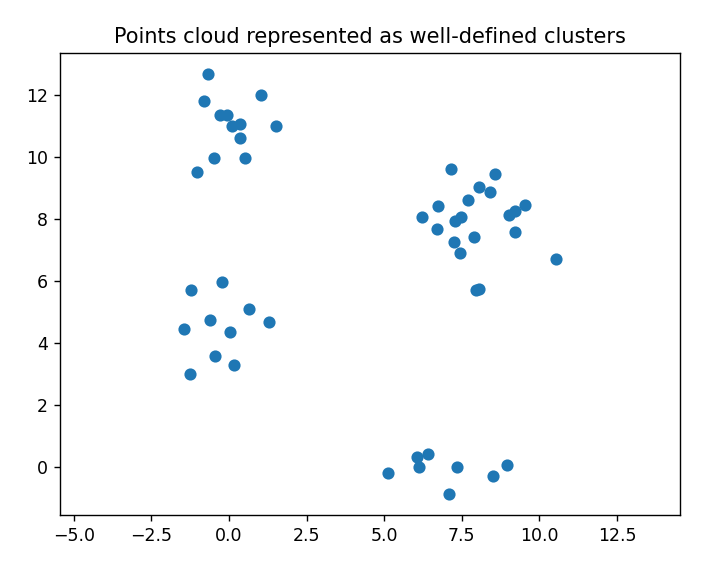}  
        \caption{Well-defined clusters points cloud}
        \label{fig:cluster}
    \end{subfigure}
    \hfill
    \begin{subfigure}{0.32\textwidth}
        \centering
        \includegraphics[width=\linewidth]{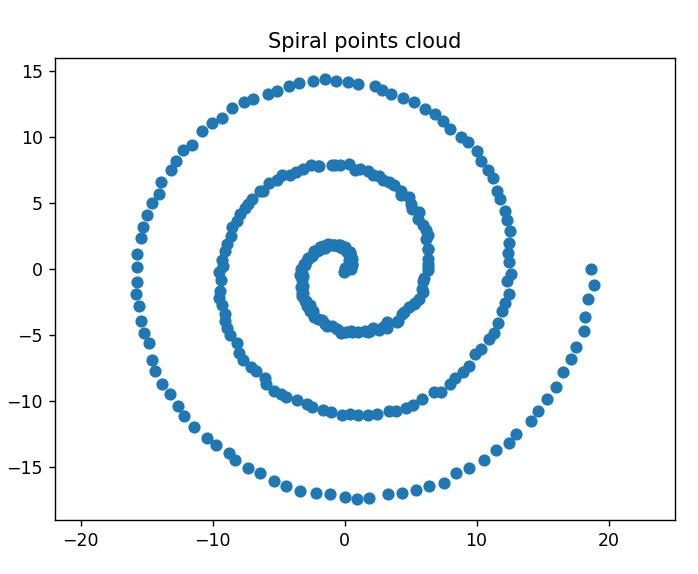}  
        \caption{Spiral points cloud}
        \label{fig:spiral}
    \end{subfigure}
    
    \caption{Visualization of three different types of point clouds.}
    \label{fig:point cloud}
\end{figure}

To analyze point clouds, we need a way to measure how "close" or "far apart" the points are from each other. For this, we rely on the idea of a \textbf{metric space}, where the concept of distance is clearly defined. By knowing how to calculate distances between points, we can begin to detect patterns, shapes, and connections within the cloud.
In the next section, we will introduce \textbf{metric spaces} more formally, to understand how distances allow us to explore and analyze these point clouds effectively.

\subsubsection{Metric Spaces}


To better understand how we measure distances between points in a point cloud, we need to define what we refer to as a \textbf{metric space}.
At its core, a metric space is simply a set of points where a clear rule—called a \textit{metric}—tells us how far apart any two points are. This notion of distance is essential for analyzing data, especially when we want to study shapes, patterns, or connections within a dataset.

Here is a simple definition:

\begin{definition}[Metric Space]
A \textbf{metric space} is a pair $(M, d)$ where:
\begin{itemize}
    \item $M$ is a set of points,
    \item $d: M \times M \rightarrow \mathbb{R}$ is a function that assigns a distance between any two points in $M$.
\end{itemize}
The function $d$ is called a \textbf{metric} if it satisfies the following four properties for all points $x, y, z \in M$:
\begin{enumerate}
    \item \textbf{Identity:} $d(x, x) = 0$ \hfill (The distance from a point to itself is zero.)
    \item \textbf{Positivity:} $d(x, y) > 0$ if $x \neq y$ \hfill (Distances are always positive between distinct points.)
    \item \textbf{Symmetry:} $d(x, y) = d(y, x)$ \hfill (The distance from $x$ to $y$ is the same as from $y$ to $x$.)
    \item \textbf{Triangle Inequality:} $d(x, z) \leq d(x, y) + d(y, z)$ \hfill (The direct path is always the shortest.)
\end{enumerate}
If these conditions are satisfied, we say that $(M, d)$ is a metric space.
\end{definition}

\vspace{0.3cm}

To make this more concrete, let’s look at a familiar example.

\begin{example}[Metric on the Real Line]
Consider the set of real numbers $M = \mathbb{R}$. We can define the distance between any two points $x$ and $y$ as:
\[
d(x, y) = |x - y|,
\]
where $|\cdot|$ denotes the absolute value.

This is the usual way we measure distance on a number line. Let’s briefly check that this function satisfies the four properties of a metric:

\begin{itemize}
    \item \textbf{Identity:} For any $x \in \mathbb{R}$, we have $d(x, x) = |x - x| = 0$.
    \item \textbf{Positivity:} If $x \neq y$, then $|x - y| > 0$.
    \item \textbf{Symmetry:} $|x - y| = |y - x|$.
    \item \textbf{Triangle Inequality:} For any $x, y, z \in \mathbb{R}$,
    \[
    |x - z| \leq |x - y| + |y - z|.
    \]
\end{itemize}
Since all four conditions are satisfied, $(\mathbb{R}, |\cdot|)$ is indeed a metric space.
\end{example}

\begin{figure}[H]
        \centering
    \includegraphics[width=0.8\linewidth]{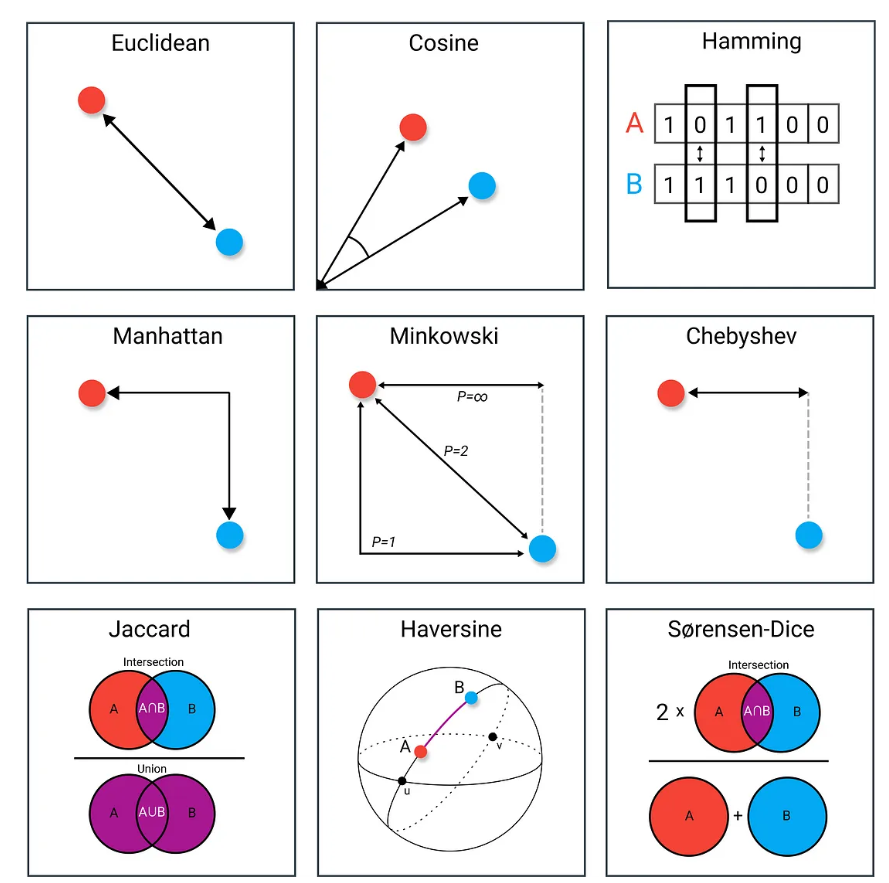}
        \caption{Some metrics commonly used in Machine Learning.  source: \href{https://medium.com/data-science/9-distance-measures-in-data-science-918109d069fa}{Maarten Grootendorst}}
        \label{metric space}
    \end{figure}

\begin{example}[Metric]
  Figure \ref{metric space} illustrates various distance measures. In our case, we specifically use the Euclidean distance.  
\end{example}


Metric spaces allow us to bring structure to sets like point clouds by giving us a precise way to talk about distances. This will be essential when we start building connections between points to uncover the hidden shape of data in Topological Data Analysis.

Now that we understand how to represent data as a \textbf{point cloud} and how to measure distances between points using \textbf{metric spaces}, the next step is to introduce some structure to these scattered points. 
While point clouds give us a raw picture of where data points are located, they don't tell us how these points are connected or how they form meaningful shapes. To capture this hidden structure, we need mathematical tools that allow us to build shapes from points in a systematic way.
This is where the concepts of \textbf{simplicial complexes} and \textbf{Homology} come into play.

\section{Simplicial Complexes and Homology}
In this section, we introduce the notions of \textbf{simplicial complexes} and \textbf{Homology}, which form the fundamental building blocks in Topological Data Analysis. Our primary focus will be understanding what simplicial complexes are and how they can be constructed from point clouds to reveal the underlying structure of data. We will then shift our focus to homology, which represents the focal point of persistent homology computation.

\subsection{Simplicial Complexes}
Simplicial complexes are mathematical structures that help us build and study the shapes or structures hidden within a set of scattered points (like a point cloud). They allow us to connect points, not just with lines, but also with higher-dimensional equivalents like triangles, tetrahedra, and beyond in a meaningful and systematic way. You can think of it as a flexible framework made up of simple pieces—points, edges, triangles, and their higher-dimensional counterparts—all glued together following specific rules.

\begin{definition}[Abstract Simplicial Complex ]\label{simplicial_complex}
Let $\mathcal{V}$ be a non-empty set of points, called \textbf{vertices}. An \textbf{Abstract Simplicial Complex} of $\mathcal{K}$ over $\mathcal{V}$ is a collection of non-empty subsets of $\mathcal{V}$ that satisfies two simple rules:
\begin{enumerate}
    \item Every single vertex $\{v\}$ with $v \in \mathcal{V}$ is in $\mathcal{K}$.
    \item If a set of vertices $\sigma$ belongs to $\mathcal{K}$, then all subsets of $\sigma$ also belongs to $\mathcal{K}$.
\end{enumerate}
\end{definition}

In simpler terms, a simplicial complex is formed by combining simplices in such a way that:
\begin{itemize}
    \item All the vertices are part of the complex.
    \item Whenever you include a shape (like a triangle) $\{v_0,v_1,v_2\}$, you must also include all its edges $\{v_0,v_1\}\{v_0,v_2\},\{v_1,v_2\}$ and points $\{v_0\},\{v_1\},\{v_2\}$.
\end{itemize}
Each building block of the complex is called a simplex. %

\subsubsection{Construction of Simplicial Complexes}

Once we understand what simplicial complexes are, the next question is: \textit{How do we build them from data?}  The construction typically begins with a \textbf{point cloud}—a set of points placed in a metric space where distances between points can be measured. In this context, each point in the cloud represents a \textbf{0-simplex}, which serves as the foundation for building higher-dimensional shapes. The next step is to create connections between these points to form higher-dimensional simplices:
\begin{itemize}
    \item When two points are "close enough", we connect them with a line segment, forming a \textbf{$1$-simplex} (an edge).
    \item If three points are all pairwise connected, we can fill in the triangle between them to form a \textbf{$2$-simplex}.
    \item Similarly, four fully connected points form a \textbf{$3$-simplex} (a tetrahedron), and so on.
\end{itemize}

But how do we decide when points are "close enough" to be connected? This is where the idea of a \textbf{threshold distance} comes in. 

\paragraph{Using distance to build connections}  
First, we calculate the distances between every pair of points in the cloud. However, these distances alone don't tell us which points should be connected. To make this decision, we introduce a threshold distance—a value that defines what we consider to be "close".If the distance between two points is less than or equal to this threshold, we draw an edge between them. As we increase the threshold distance, more points become connected, and higher-dimensional simplices begin to appear. Initially, the point cloud consists of many disconnected points, but as connections grow, the structure becomes richer and more informative.

\paragraph{Geometric representation:}  
Mathematically, a $d$-dimensional simplex is represented as the \textbf{convex hull} of its $d+1$ vertices in $\mathbb{R}^d$. This means it's the set of all weighted combinations of its vertices, where the weights are non-negative and sum to 1:
\begin{align*}
\left\lbrace \sum_{i=0}^{d} t_i v_i \ \bigg| \ t_i \geq 0, \ \sum_{i=0}^{d} t_i = 1 \right\rbrace,
\end{align*}
where $v_i$ are the vertices and $t_i$ are coefficients determining the position within the simplex.

\paragraph{\textbf{Geometric Realization of a Simplicial Complex}}  
Given a simplicial complex $K$, where $K_0$ is the set of vertices, and a mapping $\alpha: K_0 \to \mathbb{R}^N$ that places these vertices into space, the \textbf{geometric realization} of $K$, denoted by $|K|_{\alpha}$, is simply the union of all its simplices:
\begin{align*}
|K|_{\alpha} = \bigcup_{\sigma \in K} |\sigma|_{\alpha}.
\end{align*}

This process transforms an abstract set of points and connections into a concrete geometric object that approximates the shape of the data.
As we will see later, by varying the threshold distance step by step, we can create a sequence of simplicial complexes—a process called a \textbf{filtration}—which allows us to track how the shape of the data evolves across different scales.
\begin{figure}[H]
\begin{tikzpicture}[scale=2, every node/.style={circle, draw=black, fill=white, inner sep=2pt}]
\def\shift{3}

\begin{scope}[shift={(0,0,0)}]
    \node[fill=red] (A0) at (0,0,0) {};
    \node[draw=none,fill=none] at (0,-0.5,0) {0-simplex};
\end{scope}

\begin{scope}[shift={(\shift,0,0)}]
    \node[fill=blue] (A1) at (0,0,0) {};
    \node[fill=blue] (B1) at (1,0,0) {};
    \draw[thick, blue] (A1) -- (B1);
    \node[draw=none,fill=none] at (0.5,-0.5,0) {1-simplex};
\end{scope}

\begin{scope}[shift={(2*\shift,0,0)}]
    \node[fill=green] (A2) at (0,0,0) {};
    \node[fill=green] (B2) at (1,0,0) {};
    \node[fill=green] (C2) at (0.5,0.866,0) {};
    \fill[green!30,opacity=0.5] (A2.center) -- (B2.center) -- (C2.center) -- cycle;
    \draw[thick, green!70!black] (A2) -- (B2) -- (C2) -- cycle;
    \node[draw=none,fill=none] at (0.5,-0.5,0) {2-simplex};
\end{scope}

\begin{scope}[shift={(3*\shift,0,0)}]
    \node[fill=orange] (A3) at (0,0,0) {};
    \node[fill=orange] (B3) at (1,0,0) {};
    \node[fill=orange] (C3) at (0.5,0.866,0) {};
    \node[fill=orange] (D3) at (0.5,0.289,0.816) {};
    \fill[orange!30,opacity=0.5] (A3.center) -- (B3.center) -- (C3.center) -- cycle;
    \fill[orange!30,opacity=0.5] (A3.center) -- (B3.center) -- (D3.center) -- cycle;
    \fill[orange!30,opacity=0.5] (B3.center) -- (C3.center) -- (D3.center) -- cycle;
    \fill[orange!30,opacity=0.5] (A3.center) -- (C3.center) -- (D3.center) -- cycle;
    \draw[thick, orange!70!black] (A3) -- (B3) -- (C3) -- cycle;
    \draw[thick, orange!70!black] (A3) -- (D3) -- (B3) -- (D3) -- (C3) -- (D3) -- (A3);
    \node[draw=none,fill=none] at (0.5,-0.7,0) {3-simplex};
\end{scope}

\end{tikzpicture}
\end{figure}

\begin{figure}[H]
    \centering
    \begin{subfigure}{0.48\textwidth}
        \centering
        \includegraphics[width=\linewidth]{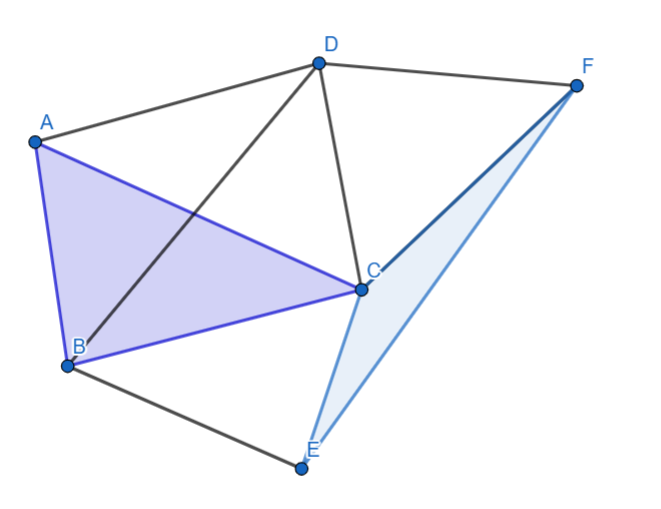} 
        \caption{2-Simplicial complex}
        \label{complex a}
    \end{subfigure}
    \hfill
    \begin{subfigure}{0.48\textwidth}
        \centering
        \includegraphics[width=\linewidth]{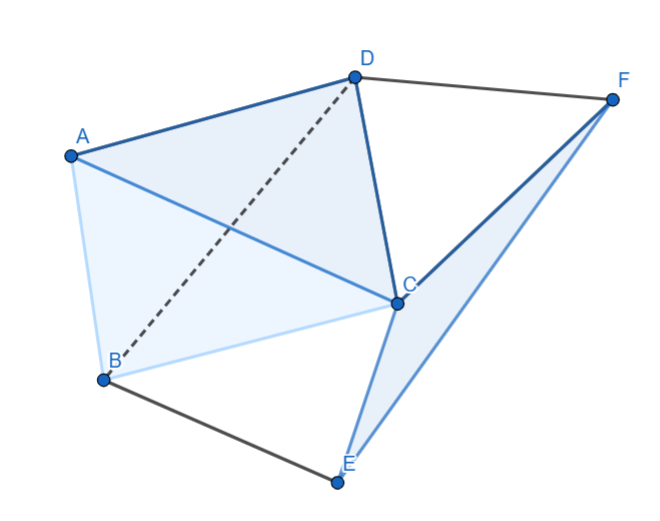}  
        \caption{3-simplicial complex}
        \label{complex b}
    \end{subfigure}
    \caption{Example of simplicial complexes}
    \label{2-3-complexes}
\end{figure}
The Figure \ref{2-3-complexes} shows a geometric realization of a $2$-simplicial complex representing the union of five 0-simplices, six 1-simplices, and two 2-simplices represented in the table below.
\begin{table}[h!]
\centering
\renewcommand{\arraystretch}{1.5} 
\setlength{\tabcolsep}{20pt} 
\fontsize{12}{15}\selectfont 
\caption{k-dimensional simplices of the complex Figure \ref{2-3-complexes}}
\begin{tabular}{|c|c|c|c|}
\hline
\textbf{k-Simplices} & $\boldsymbol{Complex ~ \ref{complex a}}$ & $\boldsymbol{Complex ~ \ref{complex b}}$  \\
\hline
\textbf{0-simplices} & \{A\},\{B\},\{C\},\{D\},\{E\},\{F\} & \{A\},\{B\},\{C\},\{D\},\{E\},\{F\}  \\
\hline
\textbf{1-simplices} & \makecell{
\{A,B\}, \{A,C\}, \{A,D\}, \{B,C\}, \\ \{B,D\}, \{B,E\},
\{C,D\}, \{C,E\},\\ \{C,F\}, \{D,F\}, \{E,F\}
} & \makecell{
\{A,B\}, \{A,C\}, \{A,D\}, \{B,C\}, \\ \{B,D\}, \{B,E\},
\{C,D\}, \{C,E\},\\ \{C,F\}, \{D,F\}, \{E,F\}
} \\
\hline
\textbf{2-simplices} & \{A,B,C\},\{C,E,F\} & \makecell{\{A,B,C\},\{A,B,D\},\{A,C,D\},\\ \{B,D,C\},\{C,E,F\}} \\
\hline
\textbf{3-simplices} & NONE & \{A,B,D,C\}  \\
\hline
\end{tabular}
\end{table}

As a simplicial complex might exhibit a high dimension, leading to high-dimensional simplices, the representation of lower-dimensional components of a $ d$-simplex is given by its faces, which are $ l$-simplices with $l<d$, where we are most interested in the $d-1$ dimensional faces.
\begin{definition}[Faces of a simplex]
   A face $\tau$ of a simplex $\sigma$ is a subset of $\sigma$  such that all vertices of $\tau$ are also vertices of $\sigma$. 
\end{definition}
A face of a simplex in a simplicial complex is a lower-dimensional simplex contained within it. Typically, we are particularly interested in the $(k-1)$-dimensional faces of a $k$-simplex, as they play a key role in defining the boundaries of the simplicial complex.
\begin{figure}[H]
    \centering
    \begin{subfigure}{0.2\textwidth}
        \centering
        \includegraphics[width=\linewidth]{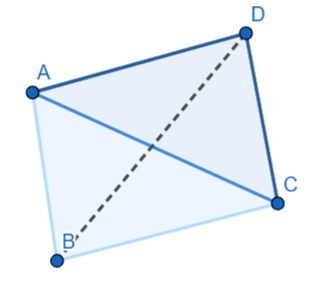}  
        \caption{Simplex}
        \label{tetra}
\end{subfigure}
    \begin{subfigure}{0.2\textwidth}
        \centering
        \includegraphics[width=\linewidth]{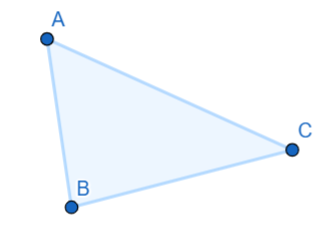} 
      \caption{Face 1}
        \label{s2}
    \end{subfigure}
    \hfill
    \begin{subfigure}{0.2\textwidth}
        \centering
        \includegraphics[width=\linewidth]{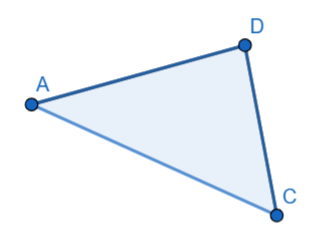}  
    \caption{Face 2}
        \label{s1}
\end{subfigure}
    \begin{subfigure}{0.18\textwidth}
        \centering
        \includegraphics[width=\linewidth]{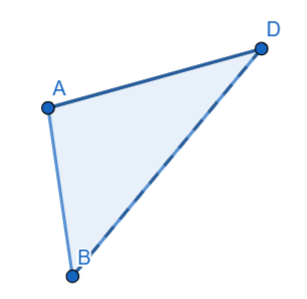} 
        \caption{Face 3}
        \label{s3}
    \end{subfigure}
    \hfill
    \begin{subfigure}{0.18\textwidth}
        \centering
        \includegraphics[width=\linewidth]{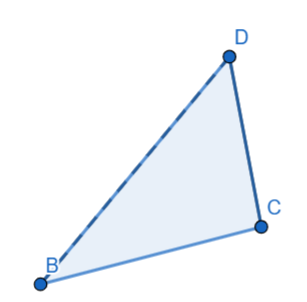}  
       \caption{Face 4}
        \label{s4}
\end{subfigure}
\caption{A Simplex with some of its faces}
\label{simplex-faces}
\end{figure}
The Figure \ref{simplex-faces} shows in \ref{tetra} a 3-simplex representing a tetrahedron with its two-dimensional faces, triangles \ref{s1}~~\ref{s2}~~\ref{s3}~~\ref{s4}.


\subsection{Homology}
Homology serves as the backbone of persistent homology, using chain complexes built from the simplices of a simplicial complex to capture the underlying structure of data across different scales. Homology reveals key topological features through the concept of chain groups and the calculation of \textbf{Betti numbers}. In this section, we will introduce the fundamental concepts of homology in a structured way, covering chain groups, boundary maps, chain complexes, homology groups and Betti Numbers.

\subsubsection{Chain Groups, Boundary Maps, and Chain Complexes.}

To understand how homology captures the shape of data, we first need to introduce a few key concepts: \textbf{chain groups}, \textbf{boundary maps}, and how they come together to form a \textbf{chain complex}.

\paragraph{\textbf{Chain Groups}}  
A chain group of a $k-$ simplicial complex \( K \) over a field $\mathbb{F}$ is a vector space with its basis formed by the elements of \( K \).\newline
For each dimension \( k \), we group together all the \( k \)-simplices of a simplicial complex \( K \). These collections form what we call the \textbf{\(k^{\text{th}}\) chain group}, denoted \( C_k(K) \). You can think of a chain group as a vector space where the basis elements are the simplices themselves, and we work over the field \( \mathbb{Z}_2 \) (where coefficients are either 0 or 1, indicating whether a simplex is included).

\begin{definition}[Chain Group]
For each dimension \( k \geq 0 \), the vector space \( C_k(K) \) over the field \( \mathbb{Z}_2 \), generated by the \( k \)-simplices of \( K \), is called the \( k^{\text{th}} \) chain group of \( K \)~\cite{nanda2021computational}.
\end{definition}

\textit{Example:}  
Consider a simplicial complex made of two triangles sharing an edge (see Figure~\ref{tetra}). The \( 2^{\text{nd}} \) chain group \( C_2(K) \) consists of all linear combinations of these two triangles. Since there are two 2-simplices, \( C_2(K) \) is isomorphic to \( \mathbb{Z}_2^2 \).

An element of a chain group, called a \textbf{\(k\)-chain}, is simply a sum of simplices:
\[
\gamma = \sum_{\sigma \in K} \gamma_{\sigma} \, \sigma, \quad \text{where } \gamma_{\sigma} \in \mathbb{Z}_2.
\]

\paragraph{\textbf{Boundary Maps}}

Boundary maps provide a way to describe how higher-dimensional simplices are connected to their faces. This is done by "breaking down" a \( k \)-simplex into its \((k-1)\)-dimensional faces.

\begin{definition}[Boundary Map]
Given the \( k^{\text{th}} \) chain group \( C_k(K) \), the \textbf{boundary map} 
\[
\delta_k^K: C_k(K) \longrightarrow C_{k-1}(K)
\]
is a linear map that sends each \( k \)-simplex to a sum of its \((k-1)\)-dimensional faces.
\end{definition}

For a simplex \( \sigma=[v_0,v_1,...,v_k] \), the boundary map is defined in general as:
\[
\delta_k^K(\sigma) = \sum_{i=0}^k (-1)^i \sigma_{-i},
\]
where \( \sigma_{-i} \) denotes the \( i^{\text{th}} \) face of \( \sigma \). However when working over the field $\mathbb{Z}/2\mathbb{Z}=\{0,1\}$, signs are supressed because $-1\cong 1 \mod 2$. This means each face appears with coefficients $0 ~or ~1$. In this case, all faces equally contribute to the formation of the complex, and a face will be cancelled out if  it appears twice as $1+1=0 \mod 2.$ Then the boundary map is simplified to:
\[
\delta_k^K(\sigma) = \sum_{i=0}^k  \sigma_{-i},
\]

Since boundary maps are linear, they act on chains by distributing over sums:
\[
\delta_k^K(\gamma) = \sum_{\sigma} \gamma_\sigma \, \delta_k^K(\sigma).
\]
\begin{example}[Boundary of a $1$-chain]\label{Boundary of a $1$-chain}

    Let's $C=[v_{0},v_{1}]+[v_{1},v_{2}]+[v_{2},v_{3}]+[v_{3},v_{4}]+[v_{4},v_{5}]+[v_{5},v_{0}]$ be a $1-chain$
    \begin{align*}
        \delta_{1}(C)&=&\delta_{1}([v_{0},v_{1}])+\delta_{1}([v_{1},v_{2}])+\delta_{1}([v_{2},v_{3}])+\delta_{1}([v_{3},v_{4}])+\delta_{1}([v_{4},v_{5}])+\delta_{1}([v_{5},v_{0}])\\
        &=&[v_{0}]+[v_{1}]+[v_{1}]+[v_{2}]+[v_{2}]+[v_{3}]+[v_{3}]+[v_{4}]+[v_{4}]+[v_{5}]+[v_{5}]+[v_{0}]\\
        &=&[v_{0}]+[v_{0}]+[v_{1}]+[v_{1}]+[v_{2}]+[v_{2}]+[v_{3}]+[v_{3}]+[v_{4}]+[v_{4}]+[v_{5}]+[v_{5}]\\
        &=&2[v_{0}]+2[v_{1}]+2[v_{2}]+2[v_{3}]+2[v_{4}]+2[v_{5}]\\
        &=&0+0+0+0+0+0=0
    \end{align*}
\end{example}
\begin{definition}[Cycle]
    Let's $C$ be a $k$-chain over a field $\mathbb{Z}/2\mathbb{Z}$. A $k-cycles$ is a $k-chain$ with \textbf{zero} boundary.  Meaning that, $C$ is a $k$-cycle if and only if $$\delta_{k}(C)=0$$.
\end{definition}
Example \ref{Boundary of a $1$-chain} demonstrates that $\delta_{1}(C) = 0$, which implies that $C$ is a $1$-cycle.

Two $1$-cycles $C$ and $C'$ are said to be homologous, or equivalent in homology, (or the same) if their sum lies in the image of the boundary map, that is, if $$C+C'\in im(\delta)$$

\begin{figure}[H]
\centering

\begin{subfigure}{0.2\textwidth}
    \centering
    \includegraphics[width=\linewidth]{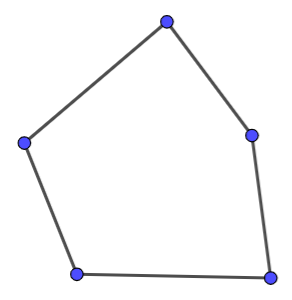}
\end{subfigure}
\hfill
\raisebox{0.8\height}{\textbf{=}}
\hfill
\begin{subfigure}{0.3\textwidth}
    \centering
    \includegraphics[width=\linewidth]{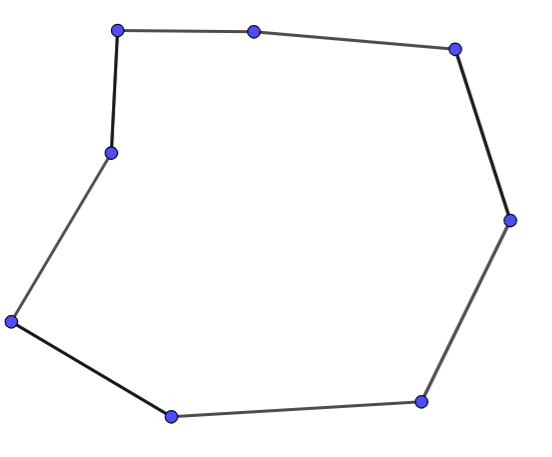}
\end{subfigure}
\hfill
\raisebox{0.8\height}{\textbf{+}}
\hfill
\begin{subfigure}{0.3\textwidth}
    \centering
    \includegraphics[width=\linewidth]{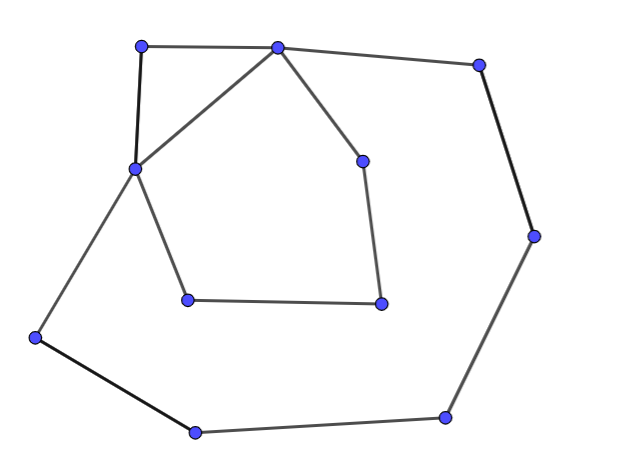}
\end{subfigure}

\end{figure}

\begin{figure}[H]
    \centering
    \[
   = \delta\left( \includegraphics[width=0.2\linewidth]{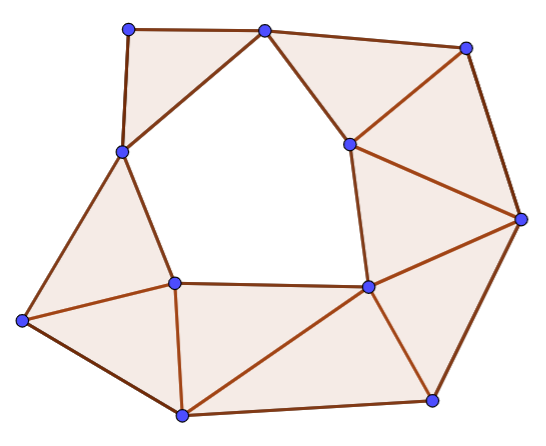} \right)
    \]
    \caption{Example of boundary of an object}
\end{figure}

\paragraph{\textbf{Chain Complexes}}
At a high level, a \textbf{chain complex} is a mathematical structure that helps us organize and analyze the different building blocks (simplices) of a simplicial complex across various dimensions. It provides a systematic way to track how these simplices connect to each other. \newline
Now that we have chain groups and boundary maps, we can link them to form a \textbf{chain complex}—a sequence of vector spaces connected by boundary maps.

\begin{definition}[Chain Complex]
A \textbf{chain complex} is a sequence of chain groups connected by linear maps (boundary maps) such that applying two boundary maps in a row always gives zero:
\[
d_{k-1} \circ d_k = 0 \quad \text{for all } k.
\]
\end{definition}

In the case of simplicial complexes, this sequence looks like:
\[
0 \longrightarrow C_k(K) \overset{\delta_k}{\longrightarrow} C_{k-1}(K) \overset{\delta_{k-1}}{\longrightarrow} \cdots \overset{\delta_1}{\longrightarrow} C_0(K) \longrightarrow 0.
\]

\begin{example}[Simplicial Chain Complex]
For a \( 3 \)-dimensional simplicial complex \( K \):
\begin{itemize}
    \item \( C_3(K) \): generated by tetrahedrons (3-simplices),
    \item \( C_2(K) \): generated by triangles (2-simplices),
    \item \( C_1(K) \): generated by edges (1-simplices),
    \item \( C_0(K) \): generated by vertices (0-simplices).
\end{itemize}
The chain complex is:
\[
0 \longrightarrow C_3(K) \overset{\delta_3}{\longrightarrow} C_2(K) \overset{\delta_2}{\longrightarrow} C_1(K) \overset{\delta_1}{\longrightarrow} C_0(K) \longrightarrow 0.
\]
\end{example}

\paragraph{\textbf{From Chains to Homology}}  
With the chain complex in place, we can now define important concepts like:
\begin{itemize}
    \item The \textbf{kernel} of a boundary map, which identifies \textbf{cycles} (closed structures with no boundary).
    \item The \textbf{image} of a boundary map, which corresponds to \textbf{boundaries} (those chains that are themselves boundaries of higher-dimensional simplices).
\end{itemize}

The difference between cycles and boundaries—formally, the kernel modulo the image—gives us the \textbf{homology groups}, which capture the essential topological features of the space. These features are quantified using \textbf{Betti numbers}, telling us how many connected components, loops, or voids exist at each dimension.

\subsubsection{Homology Groups and Betti Numbers}

In topology, one of the key challenges is understanding the essential features of a shape—those characteristics that remain unchanged even if the object is stretched, bent, or twisted. Capturing and counting these features is crucial for unveiling the intrinsic structure of the data. This is where the concepts of \textbf{homology} and \textbf{betti numbers} come into play.

\paragraph{\textbf{Homology}}
Homology provides a systematic way to detect and classify these fundamental features, known as \textbf{topological invariants}. These include:
\begin{itemize}
    \item \textbf{Connected components} (isolated pieces),
    \item \textbf{Loops} or \textbf{cycles} (like holes in a doughnut),
    \item \textbf{Voids} (empty spaces enclosed in higher dimensions),
    \item And other features in higher-dimensional spaces.
\end{itemize}

The key idea behind homology is to capture these features using algebra. This is done by constructing what are called \textbf{homology groups}.\newline

A natural question that may come to mind is \textbf{Why Homology?}\newline
Imagine deforming a rubber band into different shapes—it might stretch or twist, but as long as you don’t cut it or glue parts together, it remains a loop. Homology helps us formalize this intuition. It allows us to identify features that persist through continuous deformations, while recognizing that operations like cutting or gluing change the fundamental structure.

To study these unchanging properties, we focus on detecting "holes" of different dimensions within a space. The mathematical tool that captures this information is the \textbf{homology group}.

\begin{definition}[Homology Group]
Given a simplicial complex \( K \), the \textbf{\(k^{\text{th}}\) homology group}, denoted \( H_k(K) \), is defined as:
\[
H_k(K) = \ker(\delta_k) \, / \, \operatorname{im}(\delta_{k+1}),
\]
The homology group captures the \( k \)-dimensional features that are cycles but not boundaries, meaning they represent "holes" in dimension \( k \).
\end{definition}

The dimension of the homology group, denoted
$$\dim (H_k(K))=\dim \ker(\delta_k)-\dim \operatorname{im}(\delta_{k+1})$$
represents the number of independent $ k$-dimensional cycles which are not boundaries of a $ (k+1)$-dimensional simplices. Breaking down the formula:
\begin{itemize}
    \item \( \ker(\delta_k) \) (the kernel) counts all $k$-dimensional cycles.
    \item \( \operatorname{im}(\delta_{k+1}) \) (the image) counts all $k$-dimensional cycles that bound $(k+1)$-dimensional simplices.
    \item $\ker(\delta_k)/ \operatorname{im}(\delta_{k+1})$ counts all $k$-dimensional cycles that do not bound $(k+1)$-dimensional simplices.
\end{itemize}
This dimension is called $k^{th}$ Betti number, denoted by $\beta_k$, and it summarizes the number of $k$-dimensional holes in the space.

 \paragraph{\textbf{Betti Numbers}}  
Betti numbers count the topological features and provide a numerical summary of how many features exist in each dimension.

The \( k^{\text{th}} \) \textbf{Betti number}, denoted \( \beta_k \), counts the number of independent \( k \)-dimensional holes:
\begin{itemize}
    \item \( \beta_0 \): Number of connected components,
    \item \( \beta_1 \): Number of loops or 1-dimensional holes,
    \item \( \beta_2 \): Number of voids or enclosed surfaces,
    \item and so on for higher dimensions.
\end{itemize}

For example:
\begin{itemize}
    \item A single circle has \( \beta_0 = 1 \) (it's connected) and \( \beta_1 = 1 \) (it has one loop).
    \item A sphere has \( \beta_0 = 1 \), \( \beta_1 = 0 \) (no loops), and \( \beta_2 = 1 \) (it encloses a void).
\end{itemize}

\paragraph{Why this matters}  
By calculating homology groups and Betti numbers, we gain powerful insight into the shape of data or objects—insights that remain valid even when the object is deformed. This is the foundation of many techniques in Topological Data Analysis, where understanding the "shape" of data reveals patterns that traditional methods might miss.

In the next sections, we will see how these concepts are applied in practice, particularly through \textbf{Persistent Homology}, where we track how these topological features appear and disappear as we change the scale of observation.

\begin{table}[h!]
\centering
\renewcommand{\arraystretch}{1.5} 
\setlength{\tabcolsep}{20pt} 
\fontsize{12}{15}\selectfont 
\caption{Betti numbers for common shapes}
\begin{tabular}{|c|c|c|c|}
\hline
\textbf{Shape} & $\boldsymbol{\beta_0}$ & $\boldsymbol{\beta_1}$ & $\boldsymbol{\beta_2}$ \\
\hline
\textbf{Circle ($S^1$)} & 1 & 1 & 0 \\
\hline
\textbf{Sphere ($S^2$)} & 1 & 0 & 1 \\
\hline
\textbf{Torus} & 1 & 2 & 1 \\
\hline
\textbf{Two Disjoint Circles} & 2 & 2 & 0 \\
\hline
\end{tabular}
\end{table}

\begin{figure}[H]
    \centering
    \begin{subfigure}{0.2\textwidth}
        \centering
        \includegraphics[width=\linewidth]{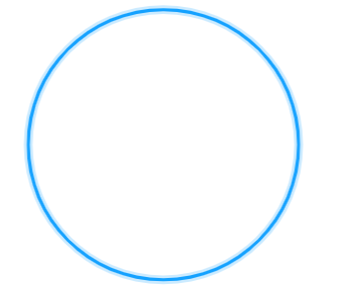} 
        \caption{Circle}
        \label{fig:circle}
    \end{subfigure}
    \hfill
    \begin{subfigure}{0.35\textwidth}
        \centering
        \includegraphics[width=\linewidth]{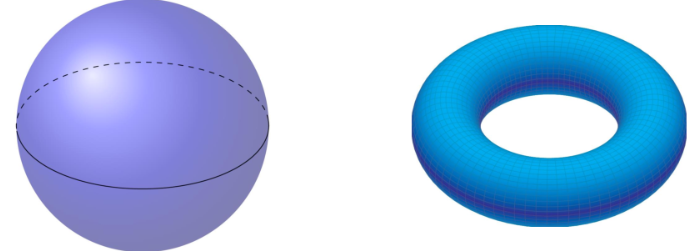}  
        \caption{Torus left and Sphere ($S^2$) right}
        \label{}
    \end{subfigure}
    \begin{subfigure}{0.42\textwidth}
        \centering
        \includegraphics[width=\linewidth]{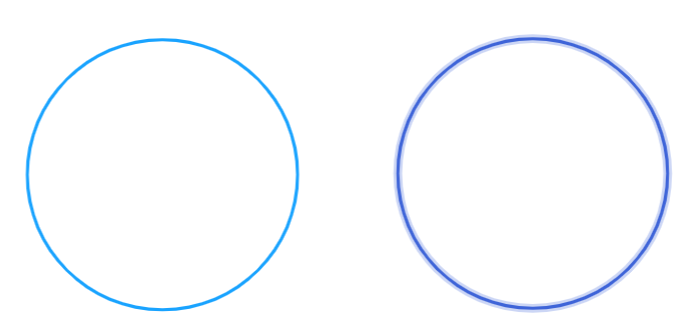}  
        \caption{Two Disjoint Circles}
        \label{small Two Disjoint Circles}
    \end{subfigure}
    \caption{Visual Representation of Common Shapes and Their Topological Features}
\end{figure}

\section{Persistent Homology}

Persistent homology (PH), a fundamental tool in \textit{Topological Data Analysis} (TDA), emerged from advancements in applied algebraic topology and computational geometry~\cite{edelsbrunner2000topological, zomorodian2005computing}. It was developed to address limitations of traditional machine learning models, particularly their inability to effectively handle incomplete or high-dimensional data, and to capture the \textit{global geometric structure} within such datasets~\cite{carlsson2009topology}. In particular, PH is very strong in dealing with noisy data, and we attribute this to the stability theorem \cite{cohen2007stability,chazal2016structure,oudot2015persistence}, which guarantees that small perturbations of the input data will lead to a small change in the resulting persistent diagram. To measure the topological differences between the persistent diagrams, distances such as the bottleneck distance and the Wasserstein distance are used.

By tracking topological features—such as connected components, loops, holes, and voids—that persist across multiple scales, PH reveals intrinsic patterns often hidden from classical statistical or geometric approaches~\cite{ghrist2008barcodes}. This robustness makes PH a valuable method in diverse fields, from biology to finance.

Persistent homology operates by constructing a nested sequence of simplicial complexes, known as a \textbf{filtration}, and computing homology at each stage. The result is a vector space whose basis consists of \textit{topological invariants}, providing a multi-scale summary of the data’s shape~\cite{edelsbrunner2010computational}. This process captures features that are stable under perturbations, offering insights into the underlying structure of complex datasets.

\subsection{Filtration and the birth/death of topological features}

This section explores how a filtration is constructed from data and how topological features emerge and disappear as the scale parameter varies.

\subsubsection{Filtration}

A \textbf{filtration} refers to a sequence of nested simplicial complexes built by gradually increasing a scale parameter $\epsilon$. For point cloud data, common filtrations include the \textit{Vietoris–Rips} and \textit{\v{C}ech} complexes~\cite{otter2017roadmap}. Formally, given a point cloud $X$, a filtration can be represented as:
\[
K_0 \subseteq K_1 \subseteq K_2 \subseteq \dots \subseteq K_n
\]
where each $K_i$ includes simplices based on proximity criteria defined by $\epsilon$. As $\epsilon$ increases, more simplices are added, enabling the detection of evolving topological features.

\subsubsection{birth/death of topological features}

During the filtration process, topological features such as connected components ($H_0$), loops ($H_1$), and voids ($H_2$) appear (\textit{birth}) and eventually disappear or merge (\textit{death}). The interval between birth and death represents the \textbf{persistence} of each feature. Long-persisting features are typically considered meaningful, while short-lived ones are often attributed to noise.

Thanks to a fundamental result in algebraic topology referred to as the structure theorem for persistence modules \cite{chazal2016structure,oudot2015persistence}, it's true that the persistent homology of a finite filtered simplicial complex can be completely characterized by a multiset of intervals. These intervals are what we visualize as barcodes, with each bar representing the lifespan of a topological feature across the filtration.

\subsection{Barcodes and Persistence diagrams}

To visualize persistent homology, two primary tools are used: \textbf{barcodes} and \textbf{persistence diagrams}~\cite{edelsbrunner2010computational}.

Barcodes represent each topological feature as a horizontal line segment, where the endpoints correspond to its birth and death scales. Longer bars indicate more persistent—and thus more significant—features.

Persistence diagrams plot features as points in a two-dimensional plane, with birth on the x-axis and death on the y-axis. Features farther from the diagonal represent higher persistence.

These visual tools not only aid interpretation but also enable integration into machine learning pipelines through vectorization techniques like \textit{persistence landscapes} and \textit{persistence images}~\cite{bubenik2015statistical, adams2017persistence}.

Persistent homology has demonstrated effectiveness in various applications, including:
\begin{enumerate}
    \item \textbf{Biology}: Protein structure and folding analysis~\cite{kasson2020topological}.
    \item \textbf{Finance}: Identifying patterns in stock market behavior~\cite{gidea2018topological}.
    \item \textbf{Neuroscience}: Studying brain network connectivity~\cite{giusti2016two}.
\end{enumerate}

Its ability to capture non-linear, multi-scale structures while remaining robust to noise makes PH an essential tool in modern data analysis.

\section{Application: DNA Structure}

To illustrate the practical power of \textbf{Persistent Homology} (PH), this section presents an example of analyzing a protein's 3D structure. Proteins, due to their complex folding patterns, are ideal candidates for topological analysis, as their shapes contain rich structural information that traditional geometric tools might overlook.

\paragraph{Step 1: Extracting Protein Data as a Point Cloud.}  
Proteins are made up of chains of amino acids, and their 3D structure is often represented by tracing the positions of specific atoms—most commonly the \textbf{alpha-carbon (C$\alpha$)} atoms along the protein backbone. These $(x, y, z)$ coordinates form a natural \textbf{point cloud}, which can then be analyzed using topological methods.

You can obtain protein 3D structures from several publicly available databases:
\begin{itemize}
    \item \href{https://alphafold.ebi.ac.uk/}{\textbf{AlphaFold Database}} -AI-predicted protein structures.
    \item \href{https://knotprot.cent.uw.edu.pl/}{\textbf{KnotProt}} - A database focused on proteins with knots and entanglements.
    \item \href{https://www.rcsb.org/}{\textbf{Protein Data Bank (PDB)}} – A comprehensive repository of experimentally determined protein structures.
\end{itemize}

Once downloaded, these structures (typically in ".pdb" format) provide the coordinates needed for topological analysis.

\paragraph{Step 2: Computing the Filtration.}  
With the point cloud in hand, the next step is to build a \textbf{filtration}—a sequence of simplicial complexes generated by connecting points based on a growing distance threshold. This process tracks how topological features (like loops and voids) appear and disappear as we change the scale.

For this, we recommend using specialized libraries such as:
\begin{itemize}
    \item \href{https://github.com/Ripser/ripser}{\textbf{Ripser}} — A fast C++/Python library for computing Vietoris–Rips persistent homology.
    \item \href{https://gudhi.inria.fr/}{\textbf{GUDHI}} — A robust Python and C++ library offering tools for simplicial complexes, filtrations, and persistent homology.
\end{itemize}

\textbf{Installation:}
\begin{verbatim}
pip install ripser
pip install gudhi
\end{verbatim}

These libraries will compute when each topological feature is \textbf{born} and when it \textbf{dies} during the filtration process. We will use the \textbf{GUDHI} library for our computations.

\paragraph{Step 3: Visualizing Topological Features.}  
The results of persistent homology are typically visualized using:
\begin{itemize}
    \item \textbf{Persistence Barcodes}: Each feature is represented by a horizontal bar, where the length indicates how long it persists.
    \item \textbf{Persistence Diagrams}: Each feature is plotted as a point in 2D space, where the x-axis represents its birth time and the y-axis its death time.
\end{itemize}

In a persistence diagram, the diagonal line $y = x$ represents features that appear and disappear almost immediately—these are usually considered \textbf{noise}. Features farther from this line are interpreted as meaningful \textbf{topological invariants}, revealing the true structural patterns within the protein.

\paragraph{Step 4: Interpreting Betti Numbers.}  
The topological features detected are quantified using \textbf{Betti numbers}:
\begin{itemize}
    \item $B_0$: Number of connected components (\(0\)-dimensional holes).
    \item $B_1$: Number of loops or cycles (\(1\)-dimensional holes).
    \item $B_2$: Number of voids (\(2\)-dimensional holes).
\end{itemize}

These numbers provide a summary of the protein's topological complexity at different dimensions.

\paragraph{Step 5: Computing Homology Generators.}  
While Betti numbers tell us \textit{how many} features exist, it's also useful to identify \textit{where} these features are located within the structure. This requires computing \textbf{homology generators}.

For this task, you can use:
\begin{itemize}
    \item \href{https://github.com/EireneML/Eirene.jl}{\textbf{Eirene.jl}} — A Julia package designed for advanced persistent homology computations.
    \item \textbf{GUDHI} — In addition to filtrations, GUDHI also supports homology generator extraction.
\end{itemize}

\textbf{Eirene.jl Installation:}
\begin{verbatim}
using Pkg
Pkg.add("Eirene")
\end{verbatim}
In our case, we will skip the last step as our focus is limited to the computation of Betti numbers.
\paragraph{Conclusion}  
By applying persistent homology to protein structures, we gain insights into the global geometric and topological properties of biomolecules—information that can be crucial in understanding protein folding, stability, and function. This method is robust to noise and dimensionality, making it a valuable tool in computational biology.

In the next section, we will walk through a complete Python implementation of this process using GUDHI and Ripser.

\subsection{Examples of persistent homology to analyze 3-1m-Supercoiled.}\label{implementation}

This paragraph will move us step by step through the computation of topological invariants of the 3-1m-supercoiled DNA molecules sequence. Since the structure of the protein sequence is known, the $3D$ structure is downloaded from the AlphaFold database\cite{varadi2022alphafold}, where we extracted the $xyz$ coordinates of carbon-alpha atoms. The coordinates will then be used to construct simplicial complexes that will be the key element of homology group computation.

\begin{enumerate}
\item \textbf{$xyz$ coordinates of carbon-alpha atoms of 3-1m-supercoiled DNA molecules.}

The analysis start with a 3-1m-supercoiled DNA molecules sequence represented as $x,y,z$ coordinates in $\mathbb{R}^{3}$
\begin{figure}[H]
    \centering
    \includegraphics[width=0.65\linewidth]{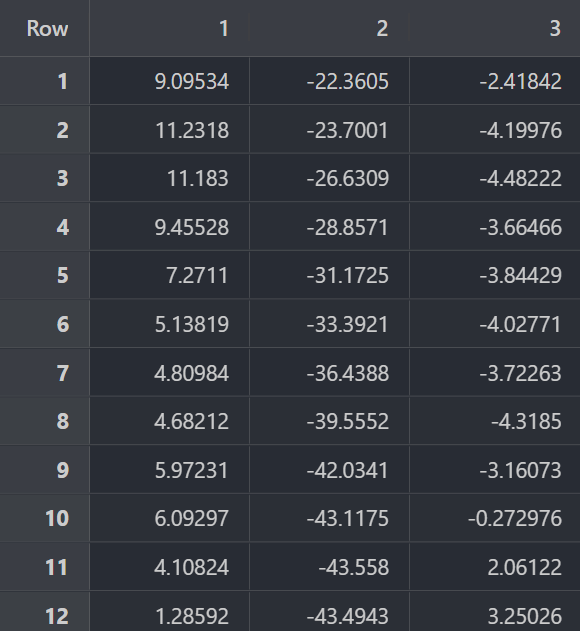}
    \caption{x,y,z coordinates of 3-1m-Supercoiled }
    \label{3-1m-Supercoiled}
\end{figure}

The figure \ref{3-1m-Supercoiled} is a table with 3 columns and 111 rows. Where the columns are respectively the $x,y,~and~z$ coordinates, and the number of rows represents the number of points used to represent the $3D$ structure of the 3-1m-Supercoiled DNA structure. Having such coordinates, the best way to give a sense of its hidden structure is to visualize them.

\item \textbf{Point Cloud}

Having the $x,y,z$ coordinates of DNA structure, we can visualize the point cloud representing those coordinates using a scatter plot as shown by the figure \ref{Point cloud of the 3-1m-Supercoiled}

\begin{figure}[H]
    \centering
    \includegraphics[width=0.7\linewidth]{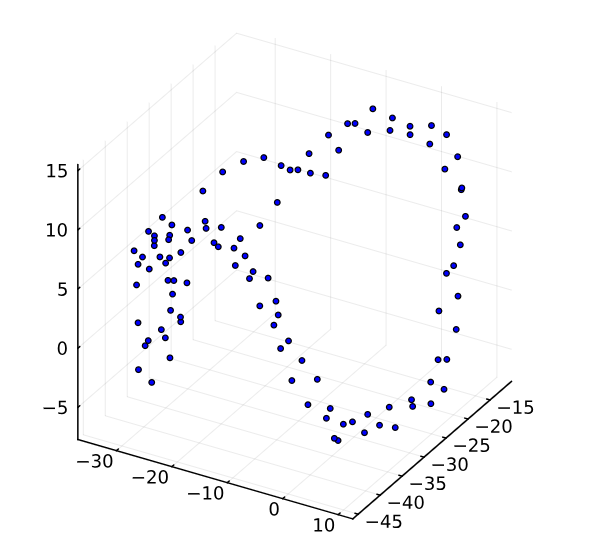}
    \caption{Point cloud of the 3-1m-Supercoiled.}
    \label{Point cloud of the 3-1m-Supercoiled}
\end{figure}

From the figure \ref{Point cloud of the 3-1m-Supercoiled}, we can see each data point with its corresponding $x,y,z$ coordinates, leading us to the notion of neighbourhood. We want to measure how data points are close or far from each other, and for that, we will start by computing pairwise distances between them, and after that, we will define our understanding of closeness by setting a threshold distance and make a comparison with the computed distances. Two points will be said to be close if the distance between them is less than or equal to the threshold distance, and then, an edge will be constructed between them, and the corresponding edge is called a 1-simplex. 

\item \textbf{Pairewise distances}

The pairwise distances are the different distances given between each pair of points. So the distance matrix obtained is a $111\times111$ matrix where rows and columns are the same, representing the different data points, the the values are the corresponding Euclidean distances obtained.

From the figure \ref{distance array}, we can see that the distance between a point and itself is zero as defined by a distance metric. 

Using the distance matrix, we can start building simplicial complexes by adding edges each time that points are close, meaning when pairwise distances are less than or equal to the threshold distance.
\begin{figure}[H]
    \centering
    \begin{subfigure}{0.5\textwidth}
        \centering
        \includegraphics[width=\linewidth]{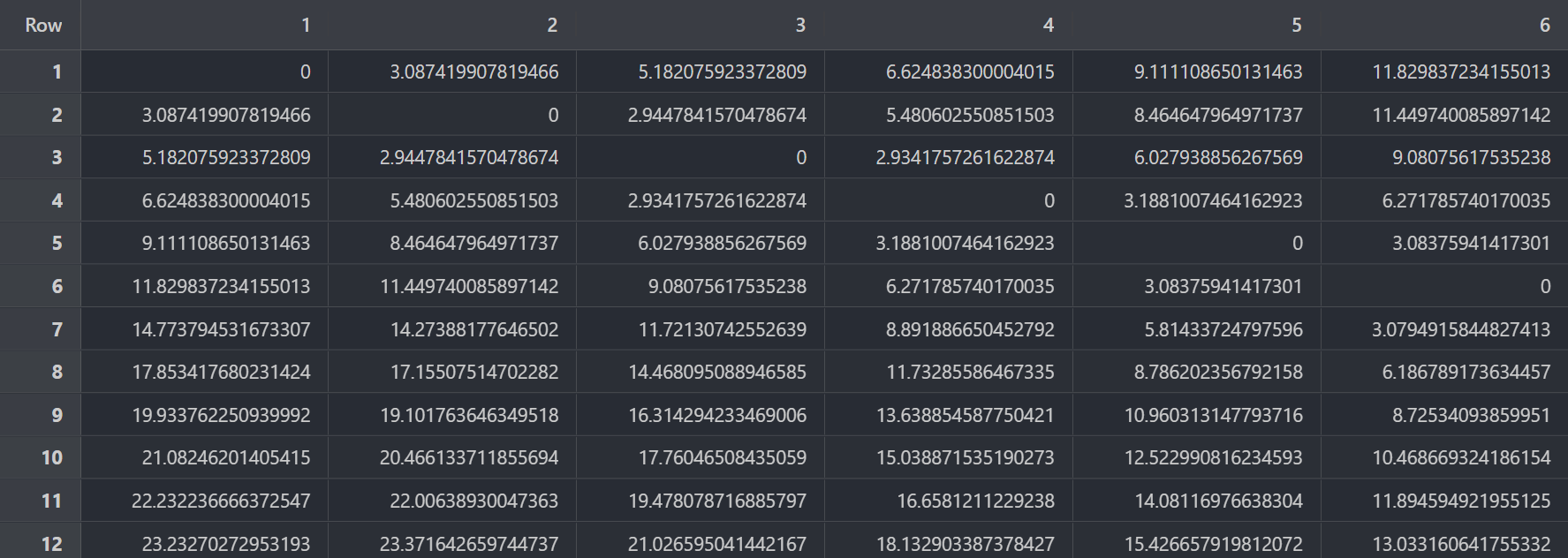} 
        \caption{pairwise distances data array}
        \label{distance array}
    \end{subfigure}
    \hfill
    \begin{subfigure}{0.45\textwidth}
        \centering
        \includegraphics[width=\linewidth]{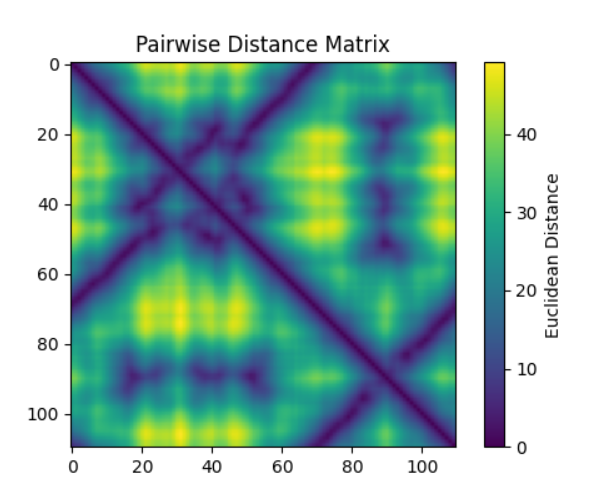}  
        \caption{Heatmap of pairwise distances}
        \label{diatanceheatmap}
    \end{subfigure}
    \caption{Pairwise distances illustrating proximity between residus}
\end{figure}

\textbf{Vietoris-Risp filtration}

The next step is to compute the filtration using Vietoris-Rips filtration as illustrated by the figure [\ref{Vietoris-Rips Filtration showing 0 and 1 simplices},\ref{Vietoris-Rips Filtration}] and showing that from a radius $r=6.5$ we start seeing a persistent hole that represents a one-dimensional hole.

\begin{figure}[H]
    \centering
    \begin{subfigure}{0.45\textwidth}
        \centering
        \includegraphics[width=\linewidth]{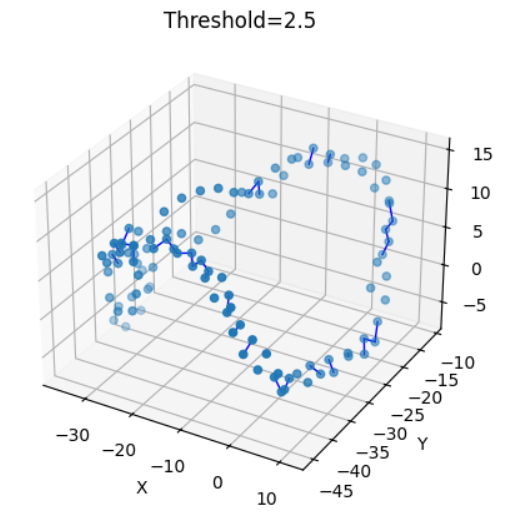} 
        \caption{}
        \label{scale 2.5}
    \end{subfigure}
    \hfill
    \begin{subfigure}{0.45\textwidth}
        \centering
        \includegraphics[width=\linewidth]{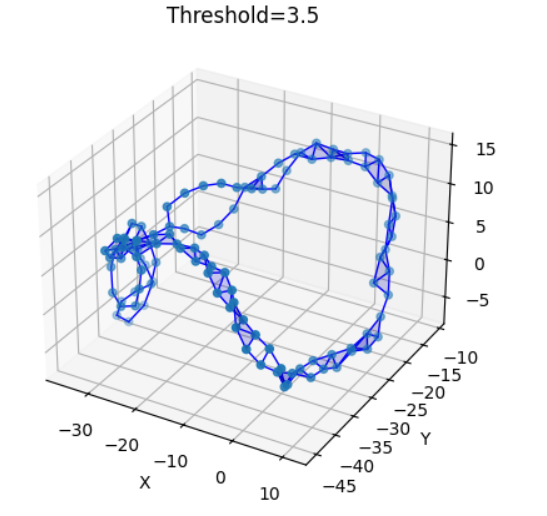}  
        \caption{}
        \label{scale 3.5}
    \end{subfigure}
    \caption{Initial Stages of Vietoris–Rips Filtration: Early Edge Formation and Emerging Loops.}
    \label{Vietoris-Rips Filtration showing 0 and 1 simplices}
\end{figure}

\begin{figure}[H]
    \centering
    \begin{subfigure}{0.45\textwidth}
        \centering
        \includegraphics[width=\linewidth]{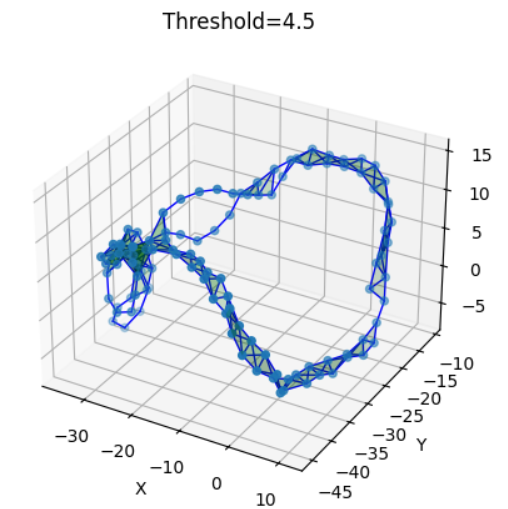} 
        \caption{}
        \label{scale 4.5}
    \end{subfigure}
    \hfill
    \begin{subfigure}{0.45\textwidth}
        \centering
        \includegraphics[width=\linewidth]{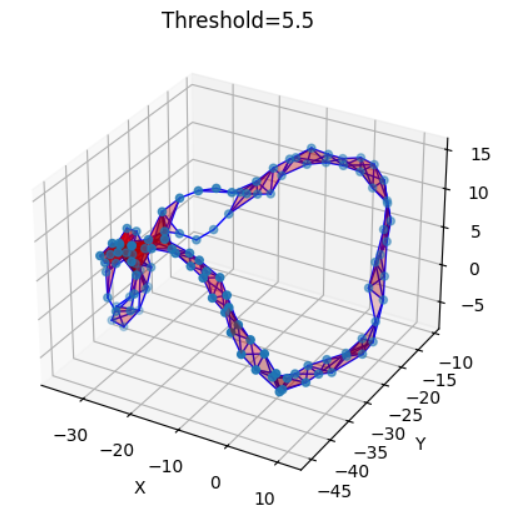}  
        \caption{}
        \label{scale 5.5}
    \end{subfigure}
    \begin{subfigure}{0.45\textwidth}
        \centering
        \includegraphics[width=\linewidth]{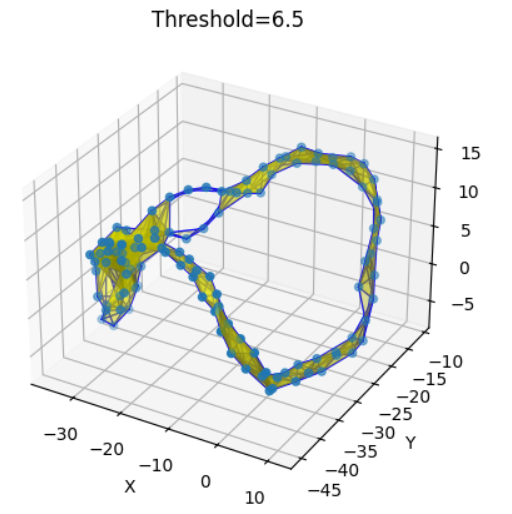} 
        \caption{}
        \label{scale 6.5}
    \end{subfigure}
    \hfill
    \begin{subfigure}{0.45\textwidth}
        \centering
        \includegraphics[width=\linewidth]{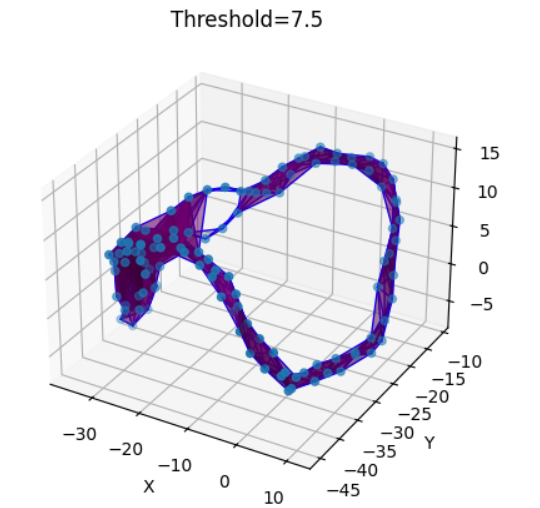}  
        \caption{}
        \label{scale 7.5}
    \end{subfigure}
    \begin{subfigure}{0.45\textwidth}
        \centering
        \includegraphics[width=\linewidth]{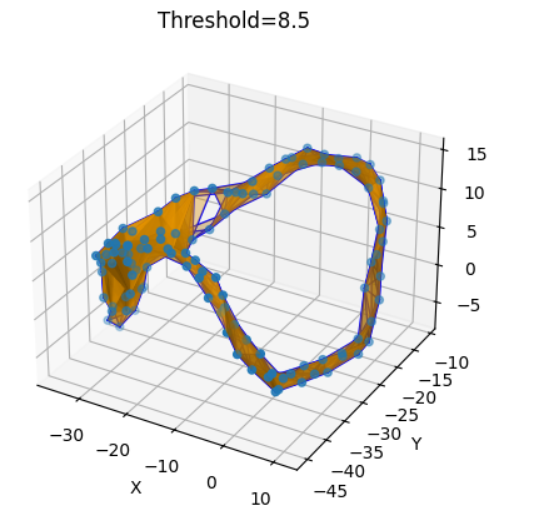} 
        \caption{}
        \label{scale 8.5}
    \end{subfigure}
    \hfill
    \begin{subfigure}{0.45\textwidth}
        \centering
        \includegraphics[width=\linewidth]{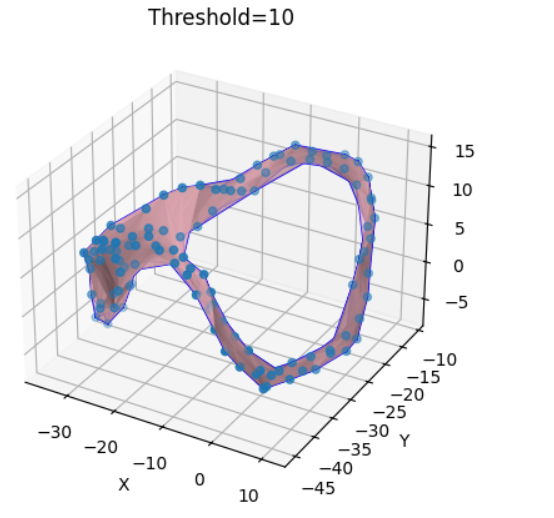}  
        \caption{}
        \label{scale 10}
    \end{subfigure}
    \caption{Vietoris–Rips Filtration Revealing the Birth and Death of 1-Dimensional Homological Features.}
    \label{Vietoris-Rips Filtration}
\end{figure}

\item \textbf{Persistent diagram}

 The persistent diagram of the filtration represented by simplicial complexes at different scales can be visualised using the persistent barcode or persistent diagram. The figure \ref{Persistent Diagram and Persistent barcode} shows on the left a persistent diagram and persistent barcode representing different homology classes, where each coordinate is defined by the birth time and the death time of each topological feature.

 \begin{figure}[H]
    \centering
    \begin{subfigure}{0.45\textwidth}
        \centering
        \includegraphics[width=\linewidth]{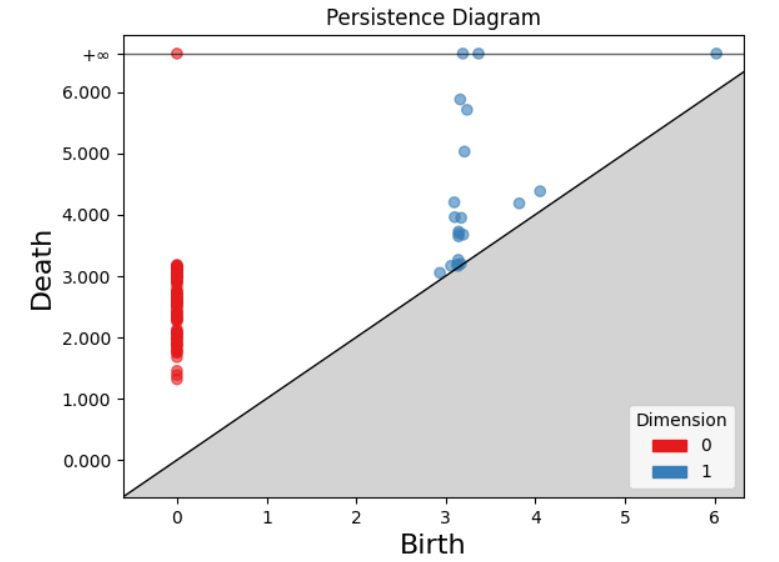} 
        \caption{Persistent Diagram}
        \label{fig:diag}
    \end{subfigure}
    \hfill
    \begin{subfigure}{0.45\textwidth}
        \centering
        \includegraphics[width=\linewidth]{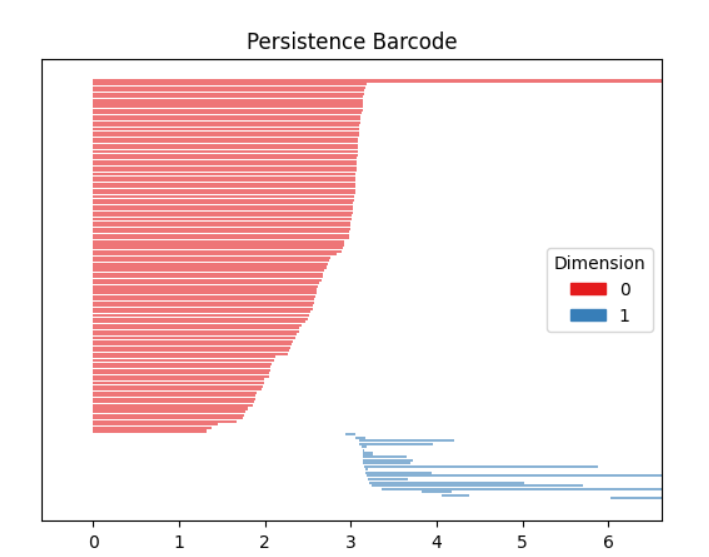}  
        \caption{Barcodes}
        \label{fig:bar}
    \end{subfigure}
    \caption{Persistence Diagram (left) and Persistence Barcode (right) showing topological features across scales computed using Vietoris-Rips Filtration. $H_0$ represents connected components, $H_1$ represents loops or cycles, and $H_2$ would represent voids (none detected here).}
    \label{Persistent Diagram and Persistent barcode}
\end{figure}

The persistence diagram shows a large number of long-lived $H_0$ features (red), indicating that many connected components persist across scales before merging. A few $H_1$ features (blue) appear but have relatively short lifespans, suggesting minor or less significant loops within the data.

Overall, the barcode confirms that the dataset is characterized by well-separated clusters or components with limited cyclic structure. This topological signature implies that machine learning algorithms focusing on connectivity or clustering could benefit from incorporating these persistent $H_0$ features, while the limited $H_1$ features may play a lesser role in downstream tasks.

 \item \textbf{Betti numbers computation}
\end{enumerate}

\begin{table}[h!]
\centering
\renewcommand{\arraystretch}{1.5} 
\setlength{\tabcolsep}{20pt} 
\fontsize{12}{15}\selectfont 
\caption{Betti Numbers for 3-1m-Supercoiled}
\begin{tabular}{|c|c|c|}
\hline
$\boldsymbol{\beta_0}$ & $\boldsymbol{\beta_1}$ & $\boldsymbol{\beta_2}$ \\
\hline
 1 & 3 & 0 \\
\hline
\end{tabular}
\end{table}

The Betti numbers presented in the above table were computed using a Python-based pipeline leveraging the \texttt{Ripser} library for topological data analysis. The numbers describe the topological structure of the 3-1m-Supercoiled DNA molecules. Specifically, $\beta_0 = 1$ indicates that the structure is composed of a single connected component. The value $\beta_1 = 3$ reveals the presence of three significant loops or cycles within the protein configuration, which may correspond to biologically relevant structural features such as twists or knots. Finally, $\beta_2=0$ confirms the absence of voids or enclosed cavities in the data, which is expected given the nature of the spatial representation of the protein.\newline
The complete Python code can be found on this \hyperlink{https://github.com/aureliejodellekemme/Persistent-Homology-from-A-to-Z.git}{GitHub repository.} 



\section{Conclusion and Future Directions}
This final section provides a systematic overview of the key points discussed and outlines potential future directions.
\subsection{Conclusion}
Persistent homology (PH) has been demonstrated to effectively complement traditional machine learning techniques by systematically uncovering the global structure of data across multiple scales. This is achieved through filtration functions, which capture topological features such as loops and voids. Furthermore, PH primarily operates on point cloud data or data that has been sampled into a point cloud, leveraging Euclidean distance to define proximity between points and construct a simplicial complex.
\subsection{Future Directions}
An interesting future direction could be exploring how persistent homology can be integrated into neural network architectures. Since both approaches involve hierarchical learning, combining them may provide deeper insights into the representations learned at each hidden layer, potentially improving interpretability and model understanding.

Another promising yet challenging future direction is the integration of Topological Data Analysis (TDA) with Graph Neural Networks (GNNs). While recent efforts have explored this intersection, existing results have not consistently demonstrated significant improvements in predictive performance over existing models \cite{taiwo2024explaining,pham2025topological}. 


\bibliographystyle{plainnat}

\section*{References}

\bibliography{paper}


\end{document}